\documentclass{amsart}
\usepackage{amssymb}
\usepackage{amsfonts}
\usepackage{latexsym}

\newtheorem{theorem}{Theorem}[section]
\newtheorem{lemma}[theorem]{Lemma}
\newtheorem{proposition}[theorem]{Proposition}
\newtheorem{corollary}[theorem]{Corollary} 
\theoremstyle{definition}

\newcommand{\Rep}{\text{Rep}}

\newcommand{\ben}{\begin{enumerate}}
\newcommand{\een}{\end{enumerate}}

\newcommand{\mC}{{\mathcal C}}

\newcommand{\CC}{{\mathbb{C}}}

\begin{document}

\title{On Vafa's theorem for tensor categories}

\author{Pavel Etingof}
\address{Department of Mathematics, Massachusetts Institute of Technology,
Cambridge, MA 02139, USA}
\email{etingof@math.mit.edu}

\maketitle

\section{Introduction}

In this note we prove two main results. 

1. In a rigid braided finite tensor category over $\CC$ 
(not necessarily semisimple), some power of the Casimir element and
some even power of the braiding is unipotent. 

2. In a (semisimple) modular category, the twists are roots of unity
dividing the algebraic integer 
$D^{5/2}$, where $D$ is the global dimension of the category
(the sum of squares of dimensions of simple objects).

Both results generalize Vafa's theorem (\cite{V,AM,BaKi}),
saying that in a modular category twists are roots of unity, and 
square of the braiding has finite order. 

In the case when the category admits a fiber functor, i.e. is a category 
of representations of a finite dimensional Hopf algebra, these results
can be found in \cite{EG1}. In fact, the method of proof of 1 and 2
is parallel to the proof of Theorems 4.3 and 4.8 in \cite{EG1},
modulo two new ingredients: categorical determinants and 
Frobenius-Perron dimensions. 

We note that statement 1 in the semisimple case was proved in \cite{Da1},
also using determinants.

At the end of the note, we discuss the notion of the quasi-exponent
of a finite rigid tensor category, which is motivated by results 1 and 2
and the paper \cite{EG3}. 

{\bf Acknowledgements.} I thank A.Davydov for very useful discussions
which inspired this note, and S.Gelaki for collaboration in \cite{EG1}, 
where the methods used here were introduced.
I am grateful to MPIM (Bonn) for hospitality. 
This research was partially supported by
the NSF grant DMS-9988796, and was done in part for the Clay Mathematics 
Institute. 

\section{Determinant of an automorphism}

Let $A$ be a 
ring, isomorphic to $\Bbb Z^N$ as an abelain group. 
 By a dimension function  
on $A$ we mean any nonzero ring homomorphism $d:A\to \Bbb C$. 
It is clear that values of $d$ are algebraic integers, since 
$d(X)$ is an eigenvalue of the matrix of left multiplication by $X$ in 
$A$ in some basis, which has integer entries.

Let $\mC$ be a tensor category over $\CC$, i.e. 
a $\CC$-linear abelian category with a biexact tensor product functor
and associativity isomorphism satisfying the pentagon identity. 
In this paper we will assume that $\mC$ is {\it finite}, i.e. equivalent,
as an abelain category, to the category of finite dimenisional 
representations of a finite dimensional $\CC$-algebra. 

Let $K_0(\mC)$ be the Grothendieck ring of $\mC$. It is a free abelian group 
of finite rank,
whose basis $S$ is the set of isomorphsim 
classes of simple objects of $\mC$, and 
multiplication is given by 
$XY=\sum_{Z\in S} [X\otimes Y:Z]Z$, where $[M:Z]$ is the multiplicity 
of occurence of $Z$ in $M$ (it is well defined by the Jordan-H\"older theorem).

Now suppose that $K_0(\mC)$ has been equipped with a dimension function 
$d$. Let $\Bbb A$ be the ring of algebraic integers.
Then to every automorphism $a:X\to X$ of an object 
$X\in \mC$ we associate its determinant, which is an element
of the $\Bbb A$-module $\Bbb A\otimes_{\Bbb Z}\Bbb C^*$, by the formula:
$\det(a)=\sum_{Z\in S}d(Z)\otimes \det(a|_{{\rm Hom}(P(Z),X)})$, where 
$P(Z)$ is the projective cover of $Z$. In the more convenient 
multiplicative notation (which we will use from now on),
$\det(a)=\prod_{Z\in S}\det(a|_{{\rm Hom}(P(Z),X)})^{d(Z)}$.
For example, if $\mC=Rep(H)$, where $H$ is a 
finite dimensional bialgebra, and $d(Z)=\dim(Z)$, then  
$\det$ takes values in 
$\Bbb Z\otimes_{\Bbb Z}\CC^*=\CC^*
\subset \Bbb A\otimes_{\Bbb Z}\CC^*$, and coincides with the usual determinant.

The properties of $\det$ are summarized in the following proposition.

\begin{proposition}\label{prope}
1) $\det(ab)=\det(a)\det(b)$. 

2) If $0\to X\to Y\to Z\to 0$ is an exact sequence respected by 
$a: Y\to Y$ then $\det(a|_Y)=\det(a|_X)\det(a|_Z)$.

3) If $a$ is a scalar then $\det(a|_X)=a^{d(X)}$.

4) $\det(1|_X\otimes a|_Y)=\det(a|_Y)^{d(X)}$.
\end{proposition}

\begin{proof}
1) is obvious. 2) follows from the fact that 
for a projective object $P$ the functor ${\rm Hom}(P,*)$ is exact.
3) follows from the fact that one has $\dim{\rm Hom}(P(Z),X)=[X:Z]$.
Now, 2) implies that it suffices to check 4) for simple 
$Y$, in which case it is clear from 3)
and multiplicativity of $d$, since $a$ is a scalar. 
\end{proof} 

{\bf Remark.} The determinant of an automorphism 
is the specialization of the class of this 
automorphism (discussed in \cite{Da1})
with respect to the dimension function $d$.

\section{Frobenius-Perron dimensions}

Recall that a finite $\Bbb Z_+$-ring
is a ring $A$ (free as a $\Bbb Z$-module) with a fixed finite basis $S$,
such that the structure constants are nonnegative.
For instance, if $\mC$ is a finite tensor category, then
$K_0(\mC)$ is a finite $\Bbb Z_+$-ring. 

We will say that $A$ is transitive if for any $X,Z\in S$ there exists 
$Y_1,Y_2\in S$ such that $XY_1$ and $Y_2X$
involve $Z$ with a nonzero coefficient. 
For example, if $\mC$ as above is a rigid monoidal category, 
then $K_0(\mC)$ is 
transitive. 

For transitive finite $\Bbb Z_+$-rings, there exists a remarkable dimension 
function $d_+$, called the Frobenius-Perron dimension. 
(For the theory of this dimension, see \cite{ENO}; 
the commutative case is discussed in \cite{FK}). 

Define the group homomorphism $d_+: R\to \Bbb C$ as follows:
for $X\in S$, let $d_+(X)$ be the maximal nonnegative eigenvalue 
of the matrix of multiplication by $X$. It exists by the Frobenius-Perron 
theorem, since this matrix has nonnegative entries. 

\begin{proposition}\label{fp}
1) $d_+$ is a ring homomorphism.

2) There exists a unique, up to scaling, 
element $R\in A_\CC$ such that 
$XR=d_+(X)R$, $X\in A$. After an appropriate 
normalization this element has positive coefficients, 
and satisfies $d_+(R)>0$ and $RY=d_+(Y)R$, $Y\in A$.  
\end{proposition}

\begin{proof}
Consider the matrix $M$ of right multiplication by $\sum_{X\in S}X$
in $A$ in the basis $S$. 
By transitivity, this matrix has strictly positive entries, 
so by the Frobenius-Perron theorem
it has a unique up to scaling eigenvector $R\in A_\CC$ with 
eigenvalue $\lambda_M$ (the maximal positive eigenvalue of $M$). 
Furthermore, this eigenvector can be normalized to have 
positive entries. 

Since $R$ is unique, it satisfies 
the equation $XR=d(X)R$ for some function $d$:
indeed, $XR$ is also an eigenvector of $M$ with eigenvalue 
$\lambda_M$, so it must be proportional to $R$. 
Furthermore, it is clear that $d$ is a dimension function. 
Since $R$ has positive entries, $d(X)=d_+(X)$ for $X\in S$. 
This implies 1). We also see that $d_+(X)>0$ for $X\in S$ (as
$R$ has positive coefficients), and hence $d_+(R)>0$.  

Now, by transitivity,
$R$ is the unique, up to scaling, solution of the system of linear equations
$XR=d_+(X)R$ (as the matrix $M'$ of left multiplication by 
$\sum_{X\in S}X$ also has positive entries). Hence, $RY=d'(Y)R$ for
some dimension function $d'$. Applying $d_+$ to both sides 
and using that $d_+(R)>0$, we find $d'=d_+$, as desired.     
\end{proof} 

The element $R$ will be called the regular object, because  
in the case of $\mC=\Rep(H)$ for a Hopf algebra $H$, 
it is indeed the class of the regular 
representation of $H$. More generally, we have the following proposition. 

\begin{proposition}\label{proj} If $\mC$ is rigid then 
$R=\sum_{X\in \mC}d_+(X)P(X)$ (up to scaling).
\end{proposition}

\begin{proof} 
We have $\sum_Xd_+(X)
\dim{\rm Hom}(P(X),Z)=d_+(Z)$ for any object $Z$. Hence, 
$$
\sum_Xd_+(X)\dim{\rm Hom}(P(X)\otimes Z,Y)=
$$
$$
\sum_Xd_+(X)\dim{\rm Hom}(P(X),Y\otimes Z^*)=
d_+(Y)d_+(Z^*)=d_+(Y)d_+(Z).
$$
 Now, 
$P(X)\otimes Z$ are projective objects. Hence, 
the formal sums $\sum_X d_+(X)P(X)\otimes Z$ and $d_+(Z)\sum_X d_+(X)P(X)$ 
are linear combinations of $P(Y)$, $Y\in \mC$, with the same coefficients. 
\end{proof} 

In the sequel, for rigid categories we will use 
the normalization of $R$ given by Proposition \ref{proj}.
In this case, 
generalizing \cite{ENO}, we call the number $d_+(R)$ the Frobenius-Perron
dimension of $\mC$. 

{\bf Remark.} Propositions \ref{fp}(2) and \ref{proj} for 
semisimple rigid categories appear in \cite{Da2}.  

\section{Generalization of Vafa's theorem}

Recall that a braiding on $\mC$ is a functorial isomorphism 
$\beta: \otimes \to \otimes^{op}$ satisfying the hexagon relation.
Recall also that in a braided rigid tensor category 
there exists an automorphism 
$z$ of the identity functor,
called the Casimir element, 
 such that 
\begin{equation}\label{casi}
z|_{X\otimes Y}=(z|_X\otimes z|_Y)
(\beta^2)^{-2}.
\end{equation}
 Namely, 
$z|_X=(u|_{^*X})^*u|_X$, where $u:Id\to **$ is the Drinfeld 
isomorphism attached to $\beta$ 
(i.e. $u|_X=
({\rm ev}_X\beta_{XX^*}\otimes 1_{X^{**}})(1_X\otimes {\rm coev}_{X^*})$, 
where ${\rm ev},{\rm coev}$ denote the evaluation and coevaluation morphisms). 

\begin{theorem}\label{genvafa} Let $\mC$ be a finite rigid 
tensor category. If $\beta$ is a braiding on $\mC$ then
for any $X,Y\in \mC$ the automorphism $\beta_{XY}^2$ is quasiunipotent, i.e. 
$(\beta_{XY}^2)^N$ is unipotent for some positive integer $N$.
Similarly, the Casimir automorphism $z_X$ is quasiunipotent.  
\end{theorem}

\begin{proof} Let us compute determinants of automorphisms 
in $\mC$ with respect to the dimension function $d_+$.  
 From hexagon relation it follows that 
$\det(\beta_{X\otimes Y,Z}^2)=
\det(\beta_{XZ}^2)^{d_+(Y)}\det(\beta_{YZ}^2)^{d_+(X)}$.
Taking $Y=R$, and using that $XR=d_+(X)R$, we derive from this that 
$\det(\beta_{XZ}^2)^{d_+(R)}=1$. 

Let us evaluate the determinant of identity (\ref{casi})
in $X\otimes R$, where $X$ is simple. We find, after cancelations: 
$$
\det(z|_X)^{d_+(R)}=\det(\beta_{XR}^2)^{-2}.
$$
Raising this to power $d_+(R)$, we get 
$\det(z|_X)^{d_+(R)^2}=1$. But $z$ is a scalar, so 
$(z|_X)^{d_+(X)d_+(R)^2}=1$. Hence, $z|_X$ is a root of unity,
and thus $z^L$ is unipotent on all objects for some $L$. 
Taking $N=2L$, we find that $(\beta_{XY}^2)^N$ is unipotent
for any $X,Y$.    
\end{proof}

\section{Refinement of Vafa's theorem for MTC}

Let now $\mC$ be a modular tensor category
(i.e, it is a semisimple ribbon category with a 
non-degenerate s-matrix). In this case, 
one has the automorphism of the identity functor $\theta$
called the twist, such that $\theta_{X\otimes Y}=
(\theta_X\otimes \theta_Y)\beta_{XY}^{-2}$
(it satisfies $\theta^2=z$). 
Vafa's theorem (see \cite{V,AM,BaKi}) claims that 
$\theta$ has finite order. 
Here we prove a refinement of this result, making an estimate on this order. 

Let $D$ be the global dimension of $\mC$, i.e. the sum of squares 
of categorical dimensions of its simple objects. It is a 
positive algebraic integer (see \cite{BaKi}). 

\begin{theorem} \label{fivehalves}
The order of $\theta$ is finite and divides $D^{5/2}$.  
\end{theorem}

\begin{proof} 
Let $d(X)$, $X\in \mC$, be the categorical dimension of $X$. 
Then $d$ is a dimension function on $K_0(\mC)$. 
Define the corresponding regular object $R$ by $R=\sum_{X\in S}d(X)X$. 
It is clear that $d(R)=D$. Let us argue as in the proof 
of Theorem \ref{genvafa}. We find that 
$\det(\beta^2_{XY})^D=1$. 
Now, taking the determinant of the equation
$\theta_{X\otimes Y}=
(\theta_X\otimes \theta_Y)\beta_{XY}^{-2}$
for $Y=R$, and raising it to power $D$, we get 
$\theta_X^{d(X)D^2}=1$. But by Lemma 1.2 in 
\cite{EG2} (see also \cite{BaKi}), 
$d(X)$ divides $D^{1/2}$. Thus, 
$\theta_X^{D^{5/2}}=1$ for any $X$, as desired. 
\end{proof}

{\bf Remark.} The power $5/2$ cannot be replaced with anything less than $2$: 
for the category of representations of $\widehat{sl}_2$ at level 1, one has 
$D=2$, while the twist of the nontrivial object is $\pm i$. 

Recall now that any modular tensor category has a central charge, 
a number $c$ defined modulo $8$. 

\begin{corollary}
Let $N$ be the number of simple objects of a modular tensor category $\mC$, 
and $D$ its global dimension. Then $cND^{5/2}/2$ is an algebraic integer.
\end{corollary}

This refines Vafa's result that $c$ is rational. 

\begin{proof}
Associated to $\mC$ are the matrices $S,T$ satisfying the equations 
$S^2=C$ (the charge conjugation matrix), and 
$(ST)^3=Ce^{\pi ic/4}$. Since $\det(C)=\pm1$, taking the determinants 
of these relations, we get 
$e^{\pi icN}=\det(T)^{12}$. But the order of $\det(T)$ 
divides $D^{5/2}$, since its eigenvalues are the twists. 
This implies the result. 
\end{proof}

\section{Quasi-exponent of a tensor category}

The above results and the work \cite{EG3} motivate 
the following definition: the {\it quasiexponent} 
${\rm qexp}(\mC)$ of 
a finite rigid tensor category $\mC$ is 
the smallest power $N$ such that $(\beta^2)^N$ is 
unipotent in the Drinfeld's center $Z(\mC)$.

\begin{proposition}\label{quas}
If $\mC=\Rep(H)$, where $H$ is a finite dimensional 
Hopf algebra, then ${\rm qexp}(\mC)$ coincides with the quasiexponent 
of $H$ defined in \cite{EG3}.
\end{proposition}

\begin{proof} By definition, 
the quasi-exponent of $H$ is the 
smallest $n$ such that $u^n$ is unipotent, where $u$ is the Drinfeld
element of the double $D(H)=H\otimes H^{*cop}$. 
On the other hand, the quasi-exponent of $\Rep(H)$ is the smallest 
$m$ such that $(R^{21}R)^m$ is unipotent, where $R$ is the universal 
R-matrix of $D(H)$. 

Recall that $\Delta(u)=(u\otimes u)(R^{21}R)^{-1}$. 
Thus $m$ is a divisor of $n$, and 
it remains to show that if $(R^{21}R)^m$ is unipotent then 
so is $u^m$. This is proved similarly to \cite{EG1}, where this is done in 
the semisimple case. 

Recall that 
we can write $u$ canonically 
as a product $u=u_{ss}u_{un}$ of commuting semisimple and 
unipotent elements. Then we have 
\begin{equation}\label{ssun}
\Delta(u_{ss}^m)\Delta(u_{un}^m)=
(u_{ss}^m\otimes u_{ss}^m)\cdot ((u_{un}^m\otimes u_{un}^m)(R^{21}R)^{-m})
\end{equation}
Note that $u_{ss}\otimes u_{ss}$ and $u_{un}\otimes u_{un}$ 
commute with $R^{21}R$, since so does $u\otimes u$. 
Hence the two sides of equation (\ref{ssun})
 are two decompositions of the same element into a product of 
commuting semisimple and unipotent parts. So these decompositions 
must coincide, which means that $u_{ss}^m$ is a grouplike element. 

By Radford's theorem \cite{R}, $u_{ss}^m=ab$, where $a,b$ are commuting 
in $D(H)$ grouplike elements of $H$ and $H^*$, respectively. 
Thus, we find that $a^{-1}b^{-1}u^m$ is unipotent, and all factors in this 
product commute. 

Now the proof is finished by the following lemma, which is a
straightforward generalization of Proposition 3.2 in \cite{EG3}:

\begin{lemma} Let $g\in H$, $\alpha\in H^*$ be grouplike elements
such that $\alpha u^mg$ is unipotent for some $m$. Then $g=\alpha=1$.  
\end{lemma}

\begin{proof} 
The proof is almost identical to the proof of Proposition 3.2
of \cite{EG3} (which treats the case $\alpha=1$). 
We have $(1- \alpha u^m g)^N=0$ for some integer $N>0$.
This is equivalent to
\begin{equation}\label{ueqn}
\sum_{k=0}^N (-1)^k \left(\begin{array}{c} N \\ k \end{array} \right)
\alpha^k u^{mk}g^{k}=0.
\end{equation}
This, in turn, is equivalent to
\begin{equation}\label{reqn}
\sum_{k=0}^N (-1)^k \left(\begin{array}{c} N \\ k \end{array} \right)R_{mk}(g^{k}\otimes
\alpha^{-k})=0,
\end{equation} 
where $R_{l}:=R(1\otimes S^2)(R)...(1\otimes S^{2l-2})(R)$. 
Indeed, equation (\ref{ueqn}) is obtained from 
equation 
(\ref{reqn}) by applying the antipode in the second component, and then 
multiplying the second component by the first component. So the equivalence 
of the two equations follows from the fact that 
the multiplication map $H^*\otimes H\to D(H)$ is a linear isomorphism.

Now, apply $1\otimes \varepsilon$ to equation (\ref{reqn}). We get 
$\sum_{k=0}^N (-1)^k
\left(\begin{array}{c} N \\ k \end{array} \right)
g^{k}=0,$ i.e., \linebreak
$(1-g)^N=0.$ However, $1-g$ is semisimple as $g$ has finite order, 
hence $g=1.$ 
Similarly, applying $\varepsilon\otimes 1$ to 
(\ref{reqn}), we get $\alpha=1$.
\end{proof}
Proposition \ref{quas} is proved.
\end{proof}

\begin{proposition}\label{propqexp}
1) ${\rm qexp}(\mC)={\rm qexp}(Z(\mC))$.

2) If $\mC$ is semisimple, ${\rm qexp}(\mC)={\rm qexp}(\mC^*)$, 
where $\mC^*$ is the dual category\footnote{About dual categories, 
see e.g. \cite{ENO}.} of $\mC$ with respect to an 
indecomposable module category.
\footnote{Thus, quasi-exponent is a Morita invariant of 
tensor categories in the sense of M\"uger.}

3) If $\mC$ is semisimple and ${\rm qexp}(\mC)=n$ the 
$(\beta^2)^n=1$. 
\end{proposition}

\begin{proof} 1) follows from the fact that $Z(Z(\mC))=Z(\mC)\otimes Z(\mC)$
as braided categories. 

2) We claim that $Z(\mC^*)$ is equivalent to 
$Z(\mC)$ with opposite braiding, as braided categories. 
Indeed, one may assume that 
$\mC={\rm Rep}(H)$, where $H$ is a weak Hopf algebra,
 and $\mC^*={\rm Rep}H^{*cop}$ (see \cite{ENO}). Then
$Z(\mC)={\rm Rep}(D(H))$, $Z(\mC^*)={\rm Rep}D(H^{*cop})$.
But $D(H)$ is isomorphic to $D(H^{*cop})$, in such a way that 
the $R$-matrix of $D(H)$ maps to the inverse of opposite $R$-matrix    
of $D(H^{*cop})$. Thus, 1) implies 2). 

To prove 3), it suffices to observe 
that identity (\ref{casi}) implies that 
$\beta^2$ is a semisimple operator on ${\rm Hom}(U,X\otimes Y)$
for any simple objects $X,Y,U$ of $Z(\mC)$ (its square is a scalar) 
\end{proof}

{\bf Remark.} In view of Proposition \ref{propqexp}(3),
in the semisimple case the quasi-exponent of $\mC$
should be called the exponent of $\mC$ (by analogy with \cite{EG1}).

\bibliographystyle{ams-alpha}

\end{document}